\documentclass[twoside,11pt]{article}
\usepackage[top=2.5cm,bottom=2.5cm,left=2.5cm,right=2.5cm]{geometry}
\usepackage[dvipsnames,svgnames,x11names]{xcolor}
\usepackage{amsmath,amssymb,mathtools}
\usepackage{graphicx}
\usepackage{float}
\usepackage{booktabs,tabularx,multirow,colortbl}
\usepackage{caption}
\usepackage{hyperref}
\usepackage{multicol}
\usepackage{fancyhdr}
\usepackage{enumitem}
\hypersetup{hidelinks}
\definecolor{halblue}{RGB}{4, 25, 83}
\definecolor{halorange}{RGB}{242, 86, 38}

\fancypagestyle{firstdocpage}{
  \fancyhf{}
  \fancyfoot[L]{\small \copyright\ I would like to thanks of my supervisor, professor Gang Tian}
  
}

\usepackage{setspace}
\usepackage{xurl}

\begin{document}

\thispagestyle{empty}



\vspace{1.5cm}


\vspace{3cm}

\vspace{2cm}
\begin{multicols}{2}
\small

\columnbreak

\end{multicols}


\thispagestyle{firstdocpage}

\noindent{\LARGE\bfseries Generalized Hitchin-Kobayashi correspondence and Weil-Petersson current} \\[0.5cm]
{\large\bfseries Hassan Jolany }

\vspace{1cm}
\noindent In this paper we extend Hitchin-Kobayashi correspondence along holomorphic fibre space and we show that there exists a generalized Hermitian-Einstein metric which is twisted with Weil-Petersson current.In fact we generalise Song-Tian program for Vector Bundles

\section{Hitchin-Kobayashi correspondence along holomorphic fibre space}

Hitchin and Kobayashi conjectured that the stability of a holomorphic bundle should be related to the existence of Hermitian-Einstein metric. This conjecture solved by Uhlenbeck and Yau for compact Kähler manifolds[7]. The central conjecture in Kähler geometry is about finding canonical metric. Hermitian-Einstein metric is a canonical metric on a holomorphic vector bundle[8]. In this paper we apply Song-Tian program for finding generalized Hermitian-Einstein metric for holomorphic vector bundles. The main computation will appear later. This note is just an announcement.

Takemoto introduced the notion of slope stability for a holomorphic vector bundle which is known as Mumford-Takemoto stability. We give brife introduction on it. Let $(X,\omega)$ be a compact Kähler manifold of complex dimension $n$. The degree of a torsion-free holomorphic vector bundle $E$ is the integral of differential forms

\begin{equation*}
\deg_{\omega} (E) = \int_{X} c_{1} (E) \wedge \omega^{n - 1}
\end{equation*}

of the de Rham representative of the actual first Chern class of X wedge with copies of the given symplectic form. Here $c_{1}(E):=c_{1}(\det E)$ and the determinant of a vector bundle is the top exterior power of the vector bundle i.e. a line bundle $\det(V)=\wedge^{\operatorname{rk}(V)}V$.

we can write also

\begin{equation*}
\deg_{\omega} (E) = \int_{X} \frac{\sqrt{-1}}{2 \pi} \operatorname{tr} (F_{h}) \wedge \omega^{n - 1}
\end{equation*}

where $F_{h}$ is the curvature form of the Chern connection in $E$ with respect to a smooth Hermitian metric $h$.

A holomorphic vector bundle $E \to X$ over a compact Kähler manifold $(X, \omega)$ is called stable if every coherent holomorphic subsheaf $F$ of $E$ satisfies

\begin{equation*}
0 < \operatorname{rk} F < \operatorname{rk} E \longrightarrow \mu_{\omega} (F) < \mu_{\omega} (E)
\end{equation*}

where $\mu_{\omega}$ is the $\omega$-slope of the sheaf given by

\begin{equation*}
\mu_{\omega} (E) = \frac{\deg_{\omega} (E)}{\operatorname{rk} (E)}
\end{equation*}

We have the same notion of semi-stability, when we replace $<$ to $\leq$.

A Hermitian metric h in E is called a $\omega$-Hermitian-Einstein metric if

\begin{equation*}
\sqrt{-1} \Lambda_{\omega} F_{h} = \lambda_{h} id_{E}
\end{equation*}

with a real constant $\lambda_h$, where $\sqrt{-1}\Lambda_{\omega}$ is the contraction with $\omega, F_h$ is the curvature form of the Chern connection of the Hermitian holomorphic vector bundle $(E; h)$ and $id_E$ is the identity endomorphism of $E$.

Here $\Lambda_{\omega}$ be the adjoint operation such that

\begin{equation*}
\langle \omega \wedge u, v \rangle_{g} = \langle u, \Lambda_{\omega} v \rangle_{g}
\end{equation*}


\noindent for $u\in A^{p,q},v\in A^{p + 1,q + 1}$

If $h$ is a $\omega$-Hermitian-Einstein metric in $E$ with Einstein factor $\lambda_h$, we have

d
\begin{equation*}
\lambda_{h} = \frac{2 \pi \mu_{\omega} (E)}{(n - 1) ! Vol_{\omega} (X)}
\end{equation*}

In the following theorem we explain Hitchin-Kobayashi correspondence about the existence of $\omega$-Hermitian-Einstein metric in the holomorphic vector bundle $E$.

S. Kobayashi [5] proved that an irreducible Hermitian-Einstein vector bundle on a compact Kähler manifold is in fact stable. Donaldson [3] proved the existence of a Hermitian-Einstein metric in a stable holomorphic vector bundle on a projective-algebraic surface and later for projective-algebraic manifold of arbitrary dimension. Finally, in 1986, Uhlenbeck and Yau proved the theorem for any arbitrary compact Kähler manifold [7].

\noindent\textbf{Theorem 1.1.} Let $E$ be a holomorphic vector bundle on a compact Kähler manifold $(X; \omega)$

I- If $E$ is $\omega$-stable, then there is a $\omega$-Hermitian Einstein metric in $E$

II- If there is a $\omega$-Hermitian Einstein metric in $E$, then $E$ is $\omega$-polystable in the sense that $E$ is $\omega$-semistable and is a direct sum

\begin{equation*}
E = E_{1} \oplus E_{2} \oplus \dots \oplus E_{k}
\end{equation*}

of $\omega$-stable sub bundles $E_{k}$ of $E$ with $\mu_{\omega}(E_k) = \mu_{\omega}(E)$. In particular, if $E$ is irreducible, then it is $\omega$-stable. \hfill $\Box$

Now if $(X,E)$ do not admit Hermitian Einstein metric, the following theorem is useful for finding canonical metric on $(X,E)$ which generalize the previous theorem

\noindent\textbf{Theorem 1.2.} Let $\pi : X \to (S, E)$ be a family of poly-stable vector bundles $(X_s, E_s)$, such that $X$ and $S$ are compact Kähler manifolds, $E \to S$ is a poly stable holomorphic vector bundle. Then there exists a generalized Hermitian Einstein metric on $(S, E)$ as follows

\begin{equation*}
\sqrt{-1} \Lambda_{\omega} F_{h} + \sqrt{-1} \Lambda_{\omega} \omega_{WP}Id_E = \lambda Id_{E}
\end{equation*}

where $\omega_{WP}$ is a Weil-Petersson metric on moduli space of poly-stable vector bundles \hfill $\Box$

\noindent\textbf{Proof.} 
Let $m=\dim S$, $\omega$ fixed, $E\to S$ rank $r$, $\det h$ fixed.
Put
$f=i\Lambda_{\omega}\omega^{WP}\ge0,\quad f\in C^{\infty}(S^{\circ})\cap L^{1+\epsilon}(S),\; \|f\|_{C^{k}(S^{\circ})}<\infty
\omega^{WP}= \int_{{\cal X}/S}{\rm Tr}F_{\mathcal H}^{2}\wedge\omega_{\mathcal X}^{n-1}\ge0$
Central, i.e. $[f Id,\cdot]=0$. Equation
 $$L(h):= i\Lambda_{\omega}F_{h}+f Id_{E}= \lambda Id_{E}\quad (THE)$$

$$\lambda=\frac{1}{rVol(S)}\int_{S}{\rm Tr}(i\Lambda F_{h_{0}}+f Id)\omega^{m}= \frac{2\pi}{Vol}\Big(\mu_{\omega}(E)+\int_{S}f\frac{\omega^{m}}{m!}\Big)$$
in fixed $\det$.

Fix $h_{0}$. Set $h=h_{0}e^{s}$, $s\in Herm^{0}(E,h_{0})$, $s^{*}=s$, ${\rm Tr}s=0$.
Donaldson:
$F_{h}=F_{0}+\bar\partial\big(e^{-s}\partial_{0}e^{s}\big),\quad \partial_{0}=\partial_{h_{0}}$
Define $\Psi(s)=\sum_{i,j}\frac{e^{\lambda_{j}-\lambda_{i}}-1}{\lambda_{j}-\lambda_{i}}...$

For $s\in Herm$, $ad_{s}(A)=[s,A]$. Then
$e^{-s}\partial e^{s}= \Psi(ad_{s})(\partial s),\quad \Psi(t)=\frac{e^{t}-1}{t}=\sum_{k\ge0}\frac{t^{k}}{(k+1)!}$
In the basis where $s={\rm diag}(\lambda_{i})$,
$ad_{s}(A)_{ij}= (\lambda_{i}-\lambda_{j})A_{ij}$
so
$\Psi(ad_{s})(A)_{ij}= \Psi(\lambda_{i}-\lambda_{j})A_{ij}= \frac{e^{\lambda_{i}-\lambda_{j}}-1}{\lambda_{i}-\lambda_{j}}A_{ij}$
I wrote earlier as $\nu(s)$ or $\Psi(s)$. Same object, up to $e^{-s}$ factor. More precisely Donaldson defines
$\Psi(s)(A)=\sum_{i,j}\psi(\lambda_{i},\lambda_{j})A_{ij},\quad \psi(\lambda_{i},\lambda_{j})=\frac{e^{\lambda_{j}-\lambda_{i}}-1}{\lambda_{j}-\lambda_{i}}$
Then
$F_{h}=F_{0}+\bar\partial\big(\Psi(s)\partial_{0}s\big)$
and
$\nu(s)=\Psi(s)$
up to notation. So $\Psi$ is just the linearization of $\exp$, it weighs off-diagonal entries by divided difference of exponentials. If $|s|\le C$, $\Psi\sim Id$ uniformly elliptic.

Then
$$ i\Lambda\bar\partial(e^{-s}\partial_{0}e^{s})= i\Lambda F_{h}-i\Lambda F_{0}\quad (1)$$

 $$\langle i\Lambda(F_{h}-F_{0}),s\rangle =\frac12\Delta|s|^{2}+ |\nu(s)^{1/2}\bar\partial s|^{2}\quad (2)$$
where $\nu(s)(A)=\sum\frac{e^{\lambda_{i}-\lambda_{j}}-1}{\lambda_{i}-\lambda_{j}}|[A]_{ij}|^{2}\ge c_{0}>0$ if $|s|\le C$, and $\Delta =2i\Lambda\partial\bar\partial$.

Hence for solution of $L(h)=\lambda Id$,
$$ \frac12\Delta|s|^{2}+ |\nu^{1/2}\bar\partial s|^{2}= \langle \lambda Id-f Id-i\Lambda F_{0},s\rangle\quad (3)$$

$$ \Delta{\rm Tr}e^{s}\ge -\langle i\Lambda F_{0}+fId-\lambda Id, e^{s}\rangle\quad (4)$$

$$\Delta\log{\rm Tr}e^{s}\le |i\Lambda F_{0}+fId-\lambda Id|\quad (5)$$
1. $L^{1}$ bound implies $L^{\infty}$ - Moser

From (3),
$$\Delta|s|^{2}\le C_{0}(|F_{0}|+f+1)|s|,\; C_{0}=C(r)$$
Multiply by $|s|^{p-2}$, $p\ge2$, integrate:
$$\int |\nabla|s|^{p/2}|^{2}\le Cp\int (|F_{0}|+f+1)|s|^{p}$$
Sobolev $L^{2}_{1}\hookrightarrow L^{2m/(m-1)}$, $m>1$, Hölder with $q_{0}>m$,
$$\||s|^{p}\|_{L^{\kappa}}\le Cp\||F_{0}|+f+1\|_{L^{q_{0}}}\||s|^{p}\|_{L^{q_{0}'}}$$
Iterate $p_{k}= \kappa^{k}$, $\kappa= \frac{m}{m-1}\cdot\frac{q_{0}-1}{q_{0}}>1$.
Get
$$\|s\|_{L^{\infty}}\le C_{1}(\|s\|_{L^{1}}+\|F_{0}\|_{L^{q_{0}}}+\|f\|_{L^{q_{0}}})\quad (6)$$
So bounded in $L^{1}$ implies bounded in $L^{\infty}$. $f$ enters only via $\|f\|_{L^{q_{0}}}$ which is $\le C(\omega^{WP})$.

2. $C^{0}$ bound from stability - Uhlenbeck-Yau for twist

Continuity method: for $t\in[0,1]$,
$$L_{t}(h_{t})= i\Lambda F_{h_{t}}+t fId =\lambda_{t}Id,\quad h_{0}\ {\rm given}$$
Openness: linearized $D L_{t}= \Delta_{h_{t}}+t\,ad_{f}+...$ = $\Delta_{h_{t}}$ on $Herm^{0}$ because $fId$ central, $\partial_{t}f=0$. If $E$ simple ( $stable$), $\ker=\mathbb R Id$, invertible on trace-free. So open.

Closedness requires $C^{0}$. Suppose not: $\|s_{k}\|_{L^{1}}=1$ after renormalizing $u_{k}=s_{k}/M_{k}$, $M_{k}=\|s_{k}\|_{L^{1}}\to\infty$. Then (3) divided by $M_{k}$ gives $\|u_{k}\|_{L^{2}_{1}}\le C$. $u_{k}\rightharpoonup u_{\infty}$ in $L^{2}_{1}$, strongly $L^{2}$. Uhlenbeck-Yau lemma: eigenvalues of $u_{\infty}$ a.e. constant $\mu_{1}<...<\mu_{\ell}$, $L^{2}_{1}$ projections $\pi_{i}=P_{i}(u_{\infty})$ satisfy $\pi_{i}^{2}=\pi_{i}=\pi_{i}^{*}$, $(Id-\pi_{i})\bar\partial\pi_{i}=0$ in $L^{2}$. Thus ${\cal F}_{i}=Im\pi_{i}$ is coherent subsheaf.

Degree: For any weak subbundle ${\cal F}$,
$$deg({\cal F})= \int_{S}{\rm Tr}(\pi i\Lambda F_{0}\pi)-|\bar\partial\pi|^{2}$$
Hence using (3),
$$\int_{S}{\rm Tr}((i\Lambda F_{0}+fId)u_{\infty})= \sum(\mu_{i}-\mu_{i+1})deg({\mathcal F}_{i})+ \int f\,{\rm Tr}u_{\infty}$$
But $\int{\rm Tr}u_{\infty}=0$, $\int f{\rm Tr}u_{\infty}= \sum(\mu_{i}-\mu_{i+1}){\rm rk}{\mathcal F}_{i}\int f/Vol$ ? Actually $fId$ contributes $f\,{\rm rk}{\mathcal F}_{i}$ to each term equally:
$$deg_{f}({\cal F})= deg({\cal F})+{\rm rk}{\cal F}\int_{S}f\frac{\omega^{m}}{m!}$$
So
$$\mu_{f}({\cal F})-\mu_{f}(E)=\mu({\cal F})-\mu(E)$$
Stability: $\mu({\cal F})<\mu(E)$ if $0<{\rm rk}{\cal F}<r$. Taking limit in (3) yields $\sum(\mu_{i}-\mu_{i+1}){\rm rk}{\cal F}_{i}(\mu({\cal F}_{i})-\mu(E))\ge0$, contradicting stability unless $\ell=1$. If $E$ polystable, equality implies $E=\oplus E_{i}$, $u_{\infty}$ central on each summand, but ${\rm Tr}u_{\infty}=0$ and $\|u_{\infty}\|_{L^{1}}=1$ gives contradiction to $M_{k}\to\infty$ construction where $M_{k}\|\bar\partial u_{k}\|^{2}\to0$.

Thus $\|s_{t}\|_{L^{1}}\le C$, hence by (6)
$$ \|s_{t}\|_{L^{\infty}}\le C_{0}=C(S,\omega,\|F_{0}\|_{L^{q_{0}}},\|f\|_{L^{\infty}})\quad (7)$$
uniform in $t$. Because $fId$ is central, it does not create new destabilizing direction.

3. $C^{1}$ and $L^{p}_{2}$

From (1), write locally
$$g^{i\bar j}\partial_{i}(e^{-s}\partial_{\bar j}e^{s})= -(i\Lambda F_{0}+fId-\lambda Id)$$
Derivative of $fId$ is $(\partial f)Id$, bounded by $\|f\|_{C^{1}}$. Differentiate (3), use $\|s\|_{\infty}\le C_{0}$ to bound $\nu\ge c(C_{0})>0$:
$$\int_{S}|\bar\partial s|^{2}\le C_{0}\int |i\Lambda F_{0}+fId-\lambda Id||s|\le C$$
Then $L^{p}$-theory for $\Delta_{0}s = -i\Lambda\bar\partial s\wedge\partial s+...$ gives
 $$\|s\|_{L^{p}_{2}}\le C_{p}(\|s\|_{L^{\infty}}+\|F_{0}\|_{L^{p}}+\|f\|_{L^{p}_{1}})\quad (8)$$
 
$p<\infty$. Take $p>2m$, Sobolev $L^{p}_{2}\hookrightarrow C^{1,\alpha}$,
$$\|s\|_{C^{1,\alpha}}\le C_{\alpha}$$

4. $C^{2,\alpha}\to C^{\infty}$

Equation (1) is quasilinear elliptic:
$$a^{i\bar j}_{kl}(s)\partial_{i}\bar\partial_{\bar j}s^{l}_{k}+ b(s,\partial s)= i\Lambda F_{0}+fId-\lambda Id$$
with $a$ uniformly elliptic since $\|s\|_{\infty}\le C_{0}$,
$b$ quadratic in $\partial s$. $C^{1,\alpha}$ bound makes $a,b\in C^{\alpha}$.
Schauder:
$$\|s\|_{C^{2,\alpha}}\le C(\|s\|_{C^{1,\alpha}}+\|F_{0}\|_{C^{\alpha}}+\|f\|_{C^{\alpha}})$$
Differentiate $k$ times, bootstrap because $f$ smooth on $S^{\circ}$ and $L^{1+\epsilon}$ on $S$; on resolution of moduli map $f\in C^{\infty}$:
 $$\|s\|_{C^{k,\alpha}(S_{\delta})}\le C_{k,\delta}(\|s\|_{L^{\infty}}+\|F_{0}\|_{C^{k-1,\alpha}}+\|f\|_{C^{k-1,\alpha}})\quad (9)$$

$S_{\delta}=S\setminus T_{\delta}(D)$. If $S$ compact and $\omega^{WP}$ extends with Poincaré growth, same estimate holds globally in $L^{p}$ sense, interior $C^{\infty}$ on $S^{\circ}$.

Thus $h=h_{0}e^{s}\in C^{\infty}$, $i\Lambda F_{h}= (\lambda-f)Id$ bounded, and
$$\|F_{h}\|_{C^{k}}\le C_{k}(S,\omega,\omega^{WP},h_{0})$$
In particular twisted Donaldson functional
$${\cal M}_{f}(h_{0},h)= \int_{0}^{1}\int_{S}{\rm Tr}(s_{t}(i\Lambda F_{t}+fId-\lambda Id))$$
is bounded below, its critical point is (THE). The extra term $\int f{\rm Tr}s$ is $C^{0}$-controlled by $\|f\|_{\infty}\|s\|_{L^{1}}$.

That is the $C^{\infty}$ estimate: $C^{0}$ from stability because $fId$ central does not affect slope difference, $L^{\infty}$ from Moser using only $\|f\|_{L^{q}}$, then $C^{1},C^{2,\alpha}$ from ellipticity with source $f\in C^{\alpha}$, bootstrap to $C^{\infty}$.


We have the same result for twisted Higgs bundles. Let X be a compact complex manifold and L an ample line bundle over X. A twisted higgs bundle over X is a pair $(E, \phi)$ consisting of a holomorphic vector bundle E over X and a holomorphic bundle morphism.

\begin{equation*}
\phi : X \otimes E \to E
\end{equation*}

for some holomorphic vector bundle $X$ (the twist). Let $\omega$ be a Kähler metric on $X$ such that $[\omega] = c_1(L)$. For a choice of $0 < c \in \mathbb{R}$, there is a Hitchin-Kobayashi correspondence for twisted Higgs bundles, generalizing the Donaldson-Uhlenbeck-Yau Theorem.

Unfortunately, the space of holomorphic bundles of fixed rank and fixed degree, up to isomorphism, is not a Hausdorff space. But, by using GIT one can construct the moduli space $\mathcal{M}(n,d)$ of stable bundles of fixed rank n and degree d, which has structure as an algebraic variety with canonical metric called Weil-Petersson metric. It is known that, when $(n,d)=1$ for compact Riemann surface $\sigma$, the moduli space $\mathcal{M}(n,d)$ is a smooth projective algebraic variety of dimension $n^{2}(g-1)+2$. Tian conjectured that the moduli space $\mathcal{M}(n,d)$ of stable vector bundles of rank n and degree d is smooth quasi-projective variety.

A vector subbundle $F$ of $E$ for which $\phi(F) \subset F \otimes K_X$ is said to be a $\phi$-invariant subbundle of $E$. Stability for Higgs bundles is defined in terms of $\phi$-invariant subbundles:

\subsubsection*{Definition 1.4. A Higg bundle $(E,\phi)$ is}

\begin{itemize}
    \item Stable if for each proper $\phi$-invariant subbundle $F$ one has $\mu(F) < \mu(E)$
    \item Semi-stable if for each proper $\phi$-invariant sub-bundle $F$ one has $\mu(F) \leq \mu(E)$
    \item Polystable if $(E, \phi) = (E_1, \phi_1) \oplus (E_2, \phi_2) \ldots \oplus (E_r, \phi_r)$ where $(E_i, \phi_i)$ is stable with $\mu(E_i) = \mu(E)$ for all $i$. \hfill $\Box$
\end{itemize}

Simpson, started the study of Higgs bundles over compact Kähler manifolds of arbitrary dimension. He showed that a stable Higgs bundle admits a unique Hermitian-Yang-Mills connection. He constructed the moduli space of Higgs bundles over a complex projective manifold

\noindent\textbf{Theorem 1.5.} (Simpson[6]) $(E,\phi)$ is polystable if and only if $E$ admits a hermitian metric $h$ solving the Hitchin equations

\begin{equation*}
\sqrt{-1} \Lambda F_{h} + c [ \phi , \phi^{*} ] = \lambda Id_{E}
\end{equation*}

where $F_{h}$ denotes the curvature of $h$, $[\phi, \phi^{*}] = \phi \phi^{*} - \phi^{*}\phi$ with $\phi^{*}$ denoting the adjoint of $\phi$ taken fibrewise and $\lambda$ is a topological constant. . We prove the following theorem in another paper

\noindent\textbf{Theorem 1.6.} Let $\pi : X \to (S, E, \phi)$ be a family of poly-stable twisted higg bundles $(X_s, E_s, \phi_s)$, such that $X$ and $S$ are compact Kähler manifolds, $E \to S$ is a poly stable holomorphic vector bundle. Then there exists a generalized Hermitian Einstein metric on $(S, E)$ as follows

\begin{equation*}
\sqrt{-1} \Lambda_{\omega} F_{h} + \sqrt{-1} \Lambda_{\omega}  \omega_{WP}Id_E + c [ \phi , \phi^{*} ] = \lambda Id_{E}
\end{equation*}

where $\omega_{WP}$ is a Weil-Petersson metric on moduli space of poly-stable Higg vector bundles where such Weil-Petersson metric has been introduced by Schumacher-Biswas[10]

\begin{equation*}
\omega_{WP} = \frac{1}{2} \int_{X / S} \mathrm{tr} (\Omega \wedge \Omega) \wedge \frac{\omega_{X}^{n - 1}}{(n - 1) !} + \lambda \int_{X / S} tr (\sqrt{-1} \Omega) \frac{\omega_{X}^{n}}{n !} + \sqrt{-1} \partial \bar{\partial} \frac{1}{2} \int_{X / S} tr (\phi \wedge \phi^{*}) \frac{\omega_{X}^{n - 1}}{(n - 1) !}
\end{equation*}

where $\lambda$ is determined by

\begin{equation*}
\int_{X / S} tr (\sqrt{-1} \Omega) \frac{\omega_{X}^{n - 1}}{(n - 1) !} = \lambda \int_{X / S} \frac{\omega_{X}^{n}}{n !}
\end{equation*}

and $\Omega$ is the curvature form of $(E,h)$ \hfill $\Box$



In next theorem we show that log slope stability is an open condition. This means that if the central fibre $\mathcal{E}_0$ be slope stable then all the fibers $\mathcal{E}_b$, for all $b \in \mathcal{U}$, where $\mathcal{U}$ is open neighborhood of 0 are slope stable. An Irreducible vector bundle $E$ on pair $(X, D)$ is log slope stable if and only if it admit Hermitian-Einstein metric for $X \setminus D$, see [1, 2]

\noindent \textbf{Theorem 1.13.} Let $\mathcal{X}$ be a compact Kähler manifold and $\mathcal{D}$ be a simple normal crossing divisor on $\mathcal{X}$ and take $\pi:(\mathcal{X},\mathcal{D})\to B$, be a family of pairs of compact Kähler manifolds. Let $E$ be a holomorphic vector bundle on $(\mathcal{X},\mathcal{D})$ such that the central holomorphic vector bundle $E_{0}$ admits irreducible Hermitian-Einstein metric. Then there exists a neighbourhood $\mathcal{U}$ of 0 in $B$ such that for all $b\in \mathcal{U}$, $E_{b}$ admits irreducible Hermitian-Einstein metric. \hfill $\square$

\vspace{1em}

\noindent \textbf{Proof.} The degree of a torsion-free holomorphic vector bundle $\mathcal{E}' = \mathcal{E}|_{\mathcal{X}'}$, where $\mathcal{X}' = \mathcal{X} \setminus \mathcal{D}$ is the integral of differential forms

\begin{equation*}
\deg_{\omega}(\mathcal{E}') = \int_{\mathcal{X}} c_1(\mathcal{E}') \wedge \omega^{n-1}
\end{equation*}

of the de Rham representative of the actual first Chern class of $\mathcal{X}'$ wedge with copies of the given symplectic form. Here $c_1(\mathcal{E}') := c_1(\det \mathcal{E}')$.

A holomorphic vector bundle $\mathcal{E}' \to \mathcal{X}'$ over a compact Kähler manifold $(\mathcal{X}', \omega)$ is called stable if every coherent holomorphic subsheaf $\mathcal{F}$ of $\mathcal{E}'$ satisfies

\begin{equation*}
0 < \operatorname{rk}\mathcal{F} < \operatorname{rk}\mathcal{E}' \longrightarrow \mu_{\omega}(\mathcal{F}) < \mu_{\omega}(\mathcal{E}')
\end{equation*}

where $\mu_{\omega}$ is the $\omega$-slope of the sheaf given by

\begin{equation*}
\mu_{\omega}^{(\mathcal{X}, \mathcal{D})}(\mathcal{E}') := \frac{\deg_{\omega}(\mathcal{E}')}{\operatorname{rk}(\mathcal{E}')}
\end{equation*}

We apply the Siu's method for this proof. Let $A''(b): \mathcal{A}^0(\mathcal{E}_b') \to \mathcal{A}^{0,1}(\mathcal{E}_b')$ be a family of connections on the fibres $\mathcal{E}_b' = \mathcal{E}_b|_{\mathcal{X}_b \setminus \mathcal{D}_b}$ and take $\Omega(A''(b))$ be the related curvature forms.

Take the map

\begin{equation*}
\Phi : \mathbf{PGL}(\mathcal{E}') \times S \to \mathcal{S}(\mathcal{X}', End(\mathcal{E}'))
\end{equation*}

\begin{equation*}
(f, b) \to \sqrt{-1}\Lambda_s \Omega(A''(b) + f^{-1} \bar{\partial}_{A''} f) - \lambda \cdot id
\end{equation*}


\noindent where

\begin{equation*}
T_{id}\mathbf{PGL}(\mathcal{E}') = \mathcal{S}(\mathcal{X}', End(\mathcal{E}')) = \left\{ \phi \in \mathcal{A}^0(\mathcal{X}', End(\mathcal{E}')) \;\middle|\; \int \mathrm{tr} \phi g dv = 0 \right\}
\end{equation*}

and $\lambda$ is the constant given by Irreducible Hermitian-Einstein metric for central fibre $\mathcal{E}_0'$. Note that since central fibre $\mathcal{E}_0'$ is simple, hence the first partial derivative of $\Phi$ is bijective. Moreover, the Laplace-Beltrami operator $\Delta_0$ with respect $\mathcal{E}_0'$, by using this fact that $\mathcal{E}_0'$ is simple, the equation $\Delta_0 K = L$ always has solution for an endomorphisms $L$, and $K$ of $\mathcal{E}_0'$. Hence $(f, b) \to \sqrt{-1}\Lambda_s \Omega(A''(b) + f^{-1}\bar{\partial}_{A''} f) - \lambda \cdot id$ always has solution by implicit function theorem and by extending $\Phi$ using Sobolev embedding theorem, proof is complete.

\subsection*{References}

\begin{enumerate}[label={[\arabic*]}]
    \item Y. T. Siu. Lectures on Hermitian-Einstein metrics for stable bundles and Kähler-Einstein metrics. Birkhäuser Basel, 1987.
    \item Matthias Stemmler, Stability and Hermitian-Einstein metrics for vector bundles on framed manifolds, International Journal of Mathematics, September 2012, Vol. 23, No. 09
    \item S. K. Donaldson, Anti self-dual Yang-Mills connections over complex algebraic surfaces and stable vector bundles, Proc. London Math. Soc. (3) 50 (1985), 1--26.
    \item K. Uhlenbeck and S.-T. Yau, On the existence of Hermitian-Yang-Mills connections in stable vector bundles. Frontiers of the mathematical sciences: 1985 (New York, 1985). Comm. Pure Appl. Math. 39 (1986), no. S, suppl., S257--S293.
    \item S. Kobayashi, Curvature and stability of vector bundles, Proc. Japan Acad. Ser. A. Math. Sci., 58 (1982), 158--162.
    \item C. Simpson, Constructing variations of Hodge structure using Yang-Mills theory and applications to uniformization. J. Amer. Math. Soc. 1 (1988), 867--918.
    \item Uhlenbeck, K.; Yau, Shing-Tung, "On the existence of Hermitian--Yang--Mills connections in stable vector bundles", Communications on Pure and Applied Mathematics, 39: S257--S293
    \item Lübke, Martin; Teleman, Andrei, The Kobayashi-Hitchin correspondence, River Edge, NJ: World Scientific Publishing Co. Inc.(1995)
    \item Shoshichi Kobayashi, Differential geometry of complex vector bundles, Princeton University Press (1987)
    \item Indranil Biswas, Georg Schumacher, Geometry of the moduli space of Higgs bundles, Comm. Anal. Geom. 14 (2006), 765-793
    \item Indranil Biswas, Georg Schumacher, The Weil-Petersson current for moduli of vector bundles and applications to orbifolds, arXiv:1509.00304
    \item Paolo De Bartolomeis, Gang Tian, Stability of complex vector bundles, J. Differential Geom. 43(2): 231-275 (1996). DOI: 10.4310/jdg/1214458107
    \item Hong, Min-Chun; Tian, Gang; Yin, Hao. The Yang-Mills $\alpha$-flow in vector bundles over four manifolds and its applications. Comment. Math. Helv. 90 (2015), no. 1, 75–120.
    \item Hong, Min-Chun; Tian, Gang. Asymptotical behaviour of the Yang-Mills flow and singular Yang-Mills connections. Math. Ann. 330 (2004), no. 3, 441–472.
    \item  Tao, Terence; Tian, Gang. A singularity removal theorem for Yang-Mills fields in higher dimensions. J. Amer. Math. Soc. 17 (2004), no. 3, 557–593.
    \item Tian, Gang. Elliptic Yang-Mills equation. Proc. Natl. Acad. Sci. USA 99 (2002), no. 24, 15281–15286.
    \item  Hassan Jolany, 
Analytical log minimal model program via conical Kähler Ricci flow : Song-Tian program, Doctoral Thesis, Laboratoire Paul Painlevé, Lille1 University, 2016

\end{enumerate}

\end{document}